# NORMAL APPROXIMATION FOR HIERARCHICAL STRUCTURES

### By Larry Goldstein

### *University of Southern California*

Given $F:[a,b]^k \to [a,b]$ and a nonconstant $X_0$ with $P(X_0 \in [a,b]) = 1$, define the hierarchical sequence of random variables $\{X_n\}_{n\geq0}$ by $X_{n+1} = F(X_{n,1}, \ldots, X_{n,k})$, where $X_{n,i}$ are i.i.d. as $X_n$. Such sequences arise from hierarchical structures which have been extensively studied in the physics literature to model, for example, the conductivity of a random medium. Under an averaging and smoothness condition on nontrivial $F$, an upper bound of the form $C\gamma^n$ for $0 < \gamma < 1$ is obtained on the Wasserstein distance between the standardized distribution of $X_n$ and the normal. The results apply, for instance, to random resistor networks and, introducing the notion of strict averaging, to hierarchical sequences generated by certain compositions. As an illustration, upper bounds on the rate of convergence to the normal are derived for the hierarchical sequence generated by the weighted diamond lattice which is shown to exhibit a full range of convergence rate behavior.

**1. Introduction.** Let $k \geq 2$ be an integer, $\mathcal{D} \subset \mathbb{R}$, $X_0$ a nonconstant random variable with $P(X_0 \in \mathcal{D}) = 1$ and $F:\mathcal{D}^k \to \mathcal{D}$ a given function. We consider the accuracy of the normal approximation for the sequence of hierarchical random variables $X_n$, where

$$(1) \qquad\qquad X_{n+1} = F(\mathbf{X}_n), \qquad n \geq 0,$$

and $\mathbf{X}_n = (X_{n,1}, \ldots, X_{n,k})^\mathsf{T}$ with $X_{n,i}$ independent, each with distribution $X_n$.

Hierarchical variables have been considered extensively in the physics literature (see [5] and the references therein), in particular to model conductivity of a random medium. The diamond lattice in particular has been considered in [3, 7]. Figure 1 shows the progression of the diamond lattice from large to small scale. At the large scale [Figure 1(a)], the system displays some conductivity along the bond between its top and bottom nodes. Inspection on a finer scale reveals the bond actually comprises four smaller bonds, each similar to Figure 1(a), connected as shown in Figure 1(b). Further inspection of each of the four bonds in Figure 1(b) reveals them to be constructed in a self-similar way from bonds at an even smaller level, giving the successive diagram Figure 1(c) and so on.

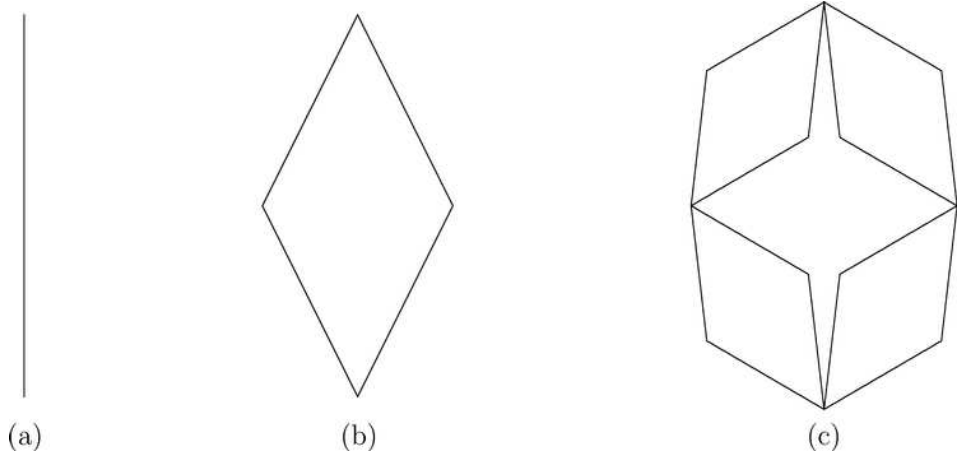

Fig. 1. *The diamond lattice.*

We assume each bond has a fixed conductivity characteristic $w \geq 0$ such that when a component with conductivity $x \geq 0$ is present along the bond the net conductivity of the bond is $wx$. For the diamond lattice as in Figure 1(b), we associate conductivities $\mathbf{w} = (w_1, w_2, w_3, w_4)^\mathsf{T}$, numbering from the top node and proceeding counterclockwise. If $\mathbf{x}_0 = (x_{0,1}, x_{0,2}, x_{0,3}, x_{0,4})^\mathsf{T}$ are the conductances of four elements each as in Figure 1(a) which are present along the bonds in Figure 1(b), then applying the resistor circuit parallel and series combination rules, the conductivity between the top and bottom nodes in Figure 1(b) is $x_1 = F(\mathbf{x}_0)$, where

$$F(\mathbf{x}) = \left( \frac{1}{w_1 x_1} + \frac{1}{w_2 x_2} \right)^{-1} + \left( \frac{1}{w_3 x_3} + \frac{1}{w_4 x_4} \right)^{-1}. \tag{2}$$

The network in Figure 1(c) is constructed from four diamond structures similar to Figure 1(b), and endowing each with the same fixed conductivity characteristics $\mathbf{w}$, with $\mathbf{x}_1 = (x_{1,1}, x_{1,2}, x_{1,3}, x_{1,4})^\mathsf{T}$ and each $x_{1,i}$ determined in the same manner as $x_1$, the conductance between the top and bottom nodes in Figure 1(c) is $x_2 = F(\mathbf{x}_1)$, and so forth.

In general, a function $F : \mathcal{D}^k \to \mathcal{D}$ and a distribution on $X_0$ such that $P(X_0 \in \mathcal{D}) = 1$ determines a sequence of distributions through $X_{n+1} =$ $F(\mathbf{X}_n)$, where $\mathbf{X}_n = (X_{n,1}, \ldots, X_{n,k})^\mathsf{T}$ with $X_{n,i}$ independent, each with distribution $X_n$. Conditions on $F$ which imply the weak law

$$(3) \qquad\qquad X_n \xrightarrow{p} c$$

have been considered by various authors. Shneiberg [8] proves that (3) holds if $\mathcal{D} = [a, b]$ and $F$ is continuous, monotonically increasing, positively homogeneous, convex and satisfies the normalization condition $F(\mathbf{1}_k) = 1$, where $\mathbf{1}_k$ is the vector of all ones in $\mathbb{R}^k$. Li and Rogers in [5] provide rather weak conditions under which (3) holds for closed $\mathcal{D} \subset (-\infty, \infty)$. See also [4, 11, 12] for an extension of the model to random $F$ and applications of hierarchical structures to computer science.

Letting $X_0$ have mean $c$ and variance $\sigma^2$, the classical central limit theorem can be set in the framework of hierarchical sequences by letting

$$(4) \qquad\qquad F(x_1, x_2) = \tfrac{1}{2}(x_1 + x_2),$$

which gives in distribution

$$X_n = \frac{X_{0,1} + \cdots + X_{0,2^n}}{2^n}.$$

Hence, $X_n \xrightarrow{p} c$, and since $X_n$ is an average of $N = 2^n$ i.i.d. variables with finite variance,

$$W_n = 2^{n/2}\left(\frac{X_n - c}{\sigma}\right) \xrightarrow{d} \mathcal{N}(0, 1).$$

Under some higher-order moment conditions one would expect a bound on the Wasserstein distance $d$ between $W_n$ and to the standard normal $\mathcal{N}$ to decay at rate $N^{-1/2}$, that is, with $\gamma = 1/\sqrt{2}$,

$$(5) \qquad\qquad d(W_n, \mathcal{N}) \le C\gamma^n.$$

The function (4), and (2) with $F(\mathbf{1}_4) = 1$, are examples of averaging functions, that is, functions $F \colon \mathcal{D}^k \to \mathcal{D}$ which satisfy the following three properties on their domain:

1. $\min_i x_i \le F(\mathbf{x}) \le \max_i x_i$.
2. $F(\mathbf{x}) \le F(\mathbf{y})$ whenever $x_i \le y_i$.
3. For all $x < y$ and for any two distinct indices $i_1 \ne i_2$, there exist $x_i \in \{x, y\}$, $i = 1, \ldots, k$, such that $x_{i_1} = x$, $x_{i_2} = y$ and $x < F(\mathbf{x}) < y$.

We note that the function $F(\mathbf{x}) = \min_i x_i$ satisfies the first two properties but not the third, and gives rise to nonnormal limiting behavior. We will call $F(\mathbf{x})$ a scaled averaging function if $F(\mathbf{x})/F(\mathbf{1}_k)$ is averaging.

Normal limits in [13] are proved for the sequences $X_n$ determined by the recursion (1) when the function $F(\mathbf{x})$ is averaging by showing that such recursions can be treated as the approximate linear recursion around the mean $c_n = EX_n$ with small perturbation $Z_n$,

$$(6) \qquad\qquad X_{n+1} = \boldsymbol{\alpha}_n \cdot \mathbf{X}_n + Z_n, \qquad n \geq 0,$$

where $\boldsymbol{\alpha}_n = F'(\mathbf{c}_n)$, $\mathbf{c}_n = (c_n, \ldots, c_n)^\mathsf{T} \in \mathbb{R}^k$ and $F'$ is the gradient of $F$. In Section 3 we prove Theorem 3.1, which gives the exponential bound (5) for the distance to the normal for sequences generated by the approximate linear recursion (6) under moment Conditions 3.1 and 3.2, which guarantee that $Z_n$ is small relative to $X_n$.

In Section 4 we prove Theorem 1.1, which shows that the normal convergence of the hierarchical sequence $X_n$ holds with the exponential bound (5) under mild conditions, and specifies $\gamma$ in an explicit range. Theorem 1.1 is proved by invoking Theorem 3.1 after showing that the required moment conditions are satisfied for averaging functions. In particular, the higher-order moment Condition 3.2 used to prove the upper bound (5) is satisfied under the same averaging assumption on $F$ used in [13] to guarantee Condition 3.1 for convergence to the normal. The condition in Theorem 1.1 that the gradient $\boldsymbol{\alpha} = F'(\mathbf{c})$ of $F$ at the limiting value $\mathbf{c}$ not be a scalar multiple of a standard basis vector rules out trivial cases such as $F(x_1, x_2) = x_1$, for which normal limits are not valid.

THEOREM 1.1. *Let $X_0$ be a nonconstant random variable with $P(X_0 \in [a, b]) = 1$ and $X_n$ given by (1) with $F : [a, b]^k \to [a, b]$, twice continuously differentiable. Suppose $F$ is averaging, or scaled averaging and homogeneous, and that $X_n \xrightarrow{p} c$, with $\boldsymbol{\alpha} = F'(\mathbf{c})$ not a scalar multiple of a standard basis vector. Then with $W_n = (X_n - c_n)/\sqrt{\mathrm{Var}(X_n)}$ and $\mathcal{N}$ a standard normal variable, for all $\gamma \in (\varphi, 1)$ there exists $C$ such that*

$$d(W_n, \mathcal{N}) \leq C\gamma^n,$$

*where*

$$(7) \qquad\qquad \varphi = \frac{\sum_{i=1}^k |\alpha_i|^3}{(\sum_{i=1}^k |\alpha_i|^2)^{3/2}},$$

*a positive number strictly less than 1. The value $\varphi$ achieves a minimum of $1/\sqrt{k}$ if and only if the components of $\boldsymbol{\alpha}$ are equal.*

At stage $n$ there are $N = k^n$ variables, so achieving the rate $\gamma^n$ for $\gamma$ to just within its minimum value $1/\sqrt{k}$ corresponds to the rate $N^{-1/2+\varepsilon}$ for every $\varepsilon > 0$. On the other hand, when $\boldsymbol{\alpha}$ is close to a standard basis vector, $\varphi$ is close to 1, and the rate $\gamma^n$ is slow. This is anticipated, as for the hierarchical sequence generated using the function, say $F(x_1, x_2) = (1 - \varepsilon)x_1 + \varepsilon x_2$ for small $\varepsilon > 0$, convergence to the normal will be slow.

In Section 5, Theorem 1.1 is applied to the hierarchical variables generated by the diamond lattice conductivity function (2). In (47) the value $\varphi$ determining the range of $\gamma$ in (5) for the rate of convergence to the normal is given as an explicit function of the weights $\mathbf{w}$; for the diamond lattice all rates $N^{-\theta}$ for $\theta \in (0, 1/2)$ are exhibited. Interestingly, there appears to be no such formula, simple or otherwise, for the limiting mean or variance of the sequence $X_n$.

We prove our results using Stein's method (see, e.g., [9]) in conjunction with the zero bias coupling of [1], derived from similar use of the size bias coupling in [2]. Let $Z$ be a mean zero, variance $\sigma^2$ normal variate and $Nh = Eh(Z/\sigma)$ for a test function $h$. Given a mean $c$, variance $\sigma^2$ random variable $X$, Stein's method, as typically applied, estimates $Eh((X - c)/\sigma) - Nh$ using the auxiliary function $f$ which is the bounded solution to

$$(8) \qquad h(w/\sigma) - Nh = \sigma^2 f'(w) - wf(w).$$

In [1] it is shown that for any mean zero variance $\sigma^2$ random variable $W$ there exists $W^*$ such that, for all absolutely continuous $f$ for which $EWf(W)$ exists,

$$(9) \qquad EWf(W) = \sigma^2 Ef'(W^*),$$

and that $W$ is normal if and only if $W \overset{d}{=} W^*$. Hence, the distance from $W$ to the normal can be expressed in a distance $d$ from $W$ to $W^*$. The variable $W^*$ is termed the $W$-zero biased distribution due to parallels with size biasing. In both size biasing and zero biasing, a sum of independent variables is biased by choosing a summand at random and replacing it with its biased version. In size biasing the variables must be nonnegative, and one is chosen with probability proportional to its expectation. In zero biasing the variables are mean zero, and one is chosen with probability proportional to its variance. The coupling construction for zero biasing just stated appears in [1] and is presented formally in Section 3; it provides the key in the proof of Lemma 2.2. To see how the zero-bias coupling is used in the Stein equation, let $f$ and $h$ be related through (8). Evaluating (8) at a mean zero, variance $\sigma^2$ variable $W$, taking expectation and using (9), we obtain

$$(10) \quad \sigma^2[Ef'(W) - Ef'(W^*)] = E[\sigma^2 f'(W) - Wf(W)] = Eh(W/\sigma) - Nh.$$

For $d$ the Wasserstein distance (also known as the Dudley, Fortet–Mourier or Kantorovich distance), Lemma 2.1 applies (10) to show the following strong connection between normal approximation and the distance between the $W$ and $W^*$ distributions as measured by $d$. With $\mathcal{N}$ a mean zero normal variable with the same variance as $W$,

$$(11) \qquad d(W, \mathcal{N}) \le 2d(W, W^*).$$

Hence, bounds on the distance between $W$ and $W^*$ can be used to bound the distance from $W$ to the normal.

We recall that, with

$$(12) \qquad \mathcal{L} = \{h \colon \mathbb{R} \to \mathbb{R} \colon |h(y) - h(x)| \leq |y - x|\},$$

the Wasserstein distance $d(Y, X)$ between variables $Y$ and $X$ on $\mathbb{R}$ is given by

$$d(Y, X) = \sup_{h \in \mathcal{L}} |E(h(Y) - h(X))|,$$

or equivalently, with

$$(13) \qquad \mathcal{F} = \{f \colon f \text{ absolutely continuous}, \; f(0) = f'(0) = 0, f' \in \mathcal{L}\}$$

we have

$$(14) \qquad d(Y, X) = \sup_{f \in \mathcal{F}} |E(f'(Y) - f'(X))|.$$

For $f \in \mathcal{F}$, certain growth restrictions are implied on $h$ of (8) for this $f$. In Theorem 3.1 these restrictions are used to compute a bound on $d(W_n, W_n^*)$, which in turn is used to bound $d(W_n, \mathcal{N})$ by (11). This argument, where $f$ is taken as given and then $h$ determined in terms of $f$ by (8), is reversed from the way Stein's method is typically applied, where $h$ is the function of interest and $f$ has only an auxiliary role as the solution of (8) for the given $h$.

For the application of Theorem 1.1, it is necessary to verify the function $F(\mathbf{x})$ in (1) is averaging. Proposition 3 of [13] shows that the effective conductance of a resistor network is an averaging function of the conductances of its individual components. Theorem 1.2, proved in Section 6, provides an additional source of averaging functions to which Theorem 1.1 may be applied by introducing the notion of strict averaging and showing that it is preserved under certain compositions.

We say $F$ is strictly averaging if strict inequality holds in property 1 when $\min_i x_i < \max_i x_i$, and in property 2 when $x_i < y_i$ for some $i$. Property 3 is the least intuitive, but is a consequence of a strict version of the first two properties; that is, a strictly averaging function is averaging: if $x < y$ and $x_{i_1} = x$, $x_{i_2} = y$, then any assignment of the values $x, y$ to the remaining coordinates gives $x < F(\mathbf{x}) < y$ by the strict form of property 1, so $F$ satisfies property 3.

THEOREM 1.2. *Let $k \geq 1$ and set $I_0 = \{1, \ldots, k\}$. Suppose subsets $I_i \subset I_0$, $i \in I_0$ satisfy $\bigcup_{i \in I_0} I_i = I_0$. For $\mathbf{x} \in \mathbb{R}^k$ and $i \in I_0$ let $\mathbf{x}_i = (x_{j_1}, \ldots, x_{j_{|I_i|}})$, where $\{j_1, \ldots, j_{|I_i|}\} = I_i$ and $j_1 < \cdots < j_{|I_i|}$. Let $(F_i \colon [0, \infty)^{|I_i|} \to [0, \infty)$ or)*

$F_i : \mathbb{R}^{|I_i|} \to \mathbb{R}$, $i = 0, \ldots, k$. If $F_0, F_1, \ldots, F_k$ are strictly averaging and $F_0$ is (positively) homogeneous, then the composition

$$F_{\mathbf{s}}(\mathbf{x}) = F_0(s_1 F_1(\mathbf{x}_1), \ldots, s_k F_k(\mathbf{x}_k))$$

is strictly averaging for any $\mathbf{s}$ for which $F_0(\mathbf{s}) = 1$ and $s_i > 0$ for all $i$. If $F_0, F_1, \ldots, F_k$ are scaled, strictly averaging and $F_0$ is (positively) homogeneous, then

$$F_{\mathbf{1}}(\mathbf{x}) = F_0(F_1(\mathbf{x}_1), \ldots, F_k(\mathbf{x}_k))$$

is a scaled strictly averaging function.

In particular, in the context of resistor networks, two components with conductances $x_1, x_2$ in parallel is equivalent to one component with conductance

$$L_1(x_1, x_2) = x_1 + x_2,$$

and in series to one component with conductance

$$L_{-1}(x_1, x_2) = (x_1^{-1} + x_2^{-1})^{-1}.$$

These parallel and series combination rules are the $p = 1$ and $p = -1$ special cases, with $w_i = 1$, of the weighted $L^p$-norm functions

$$L_p^{\mathbf{w}}(\mathbf{x}) = \left( \sum_{i=1}^{k} (w_i x_i)^p \right)^{1/p}, \qquad \mathbf{w} = (w_1, \ldots, w_k)^{\mathsf{T}}, \; w_i \in (0, \infty),$$

which are scaled, strictly averaging and positively homogeneous on $[0, \infty)^k$ for $p > 0$, and on $(0, \infty]^k$ for $p < 0$.

Though Theorem 1.2 cannot be invoked to subsume the result of [13] that every resistor network is strictly averaging in its component conductances (e.g., consider the complete graph $K_4$), now suppressing the dependence of $L_p$ on $\mathbf{w}$, since $F(\mathbf{x})$ in (2) can be represented as

$$F(\mathbf{x}) = L_1(L_{-1}(x_1, x_2), L_{-1}(x_3, x_4)),$$

Theorem 1.2 obtains to show that the diamond lattice conductivity function is a scaled, strictly averaging function on $(0, \infty)^4$ for any choice of positive weights. Moreover, Theorem 1.2 shows the same conclusion holds when the resistor parallel $L_1$ and series $L_{-1}$ combination rules in this network are replaced by, say, $L_2$ and $L_{-2}$, respectively.

**2. Zero bias and the Wasserstein distance.** The following lemma, of separate interest, shows how the zero bias coupling of $W$ upper bounds the Wasserstein distance to normality.

LEMMA 2.1. *Let $W$ be a mean zero, finite variance random variable, and let $W^*$ have the $W$-zero bias distribution. Then with $d$ the Wasserstein distance, and $\mathcal{N}$ a normal variable with the same variance as $W$,*

$$d(W, \mathcal{N}) \leq 2d(W, W^*).$$

PROOF. Since $\sigma^{-1}d(X, Y) = d(\sigma^{-1}X, \sigma^{-1}Y)$ and $\sigma^{-1}W^* = (\sigma^{-1}W)^*$, we may assume $\mathrm{Var}(W) = 1$. The dual form of the Wasserstein distance gives that

$$\inf_{(Y,X)} E|Y - X| = d(Y, X), \tag{15}$$

where the infimum, achieved for random variables on $\mathbb{R}$, is taken over all pairs $(Y, X)$ on a common space with the given marginals (see [6]). Take $W, W^*$ to achieve the infimum $d(W, W^*)$.

For a differentiable test function $h$ and $\sigma^2 = 1$, Stein [10] shows the solution $f$ of (8) is twice differentiable with $\|f''\| \leq 2\|h'\|$, where $\|\cdot\|$ represents the supremum norm. Now going from right to left in (10), applying this bound and using (15) we have

$$|Eh(W) - Nh| \leq \|f''\| E|W - W^*| \leq 2\|h'\| E|W - W^*| = 2\|h'\| d(W, W^*).$$

Functions $h \in \mathcal{L}$ of (12) are absolutely continuous with $\|h'\| \leq 1$, so taking supremum over $h \in \mathcal{L}$ on the left-hand side completes the proof. □

The following results in this section give the prototype of the argument used in Section 3 and show how the zero bias coupling can be used to obtain the exponential decay of the Wasserstein distance to the normal.

PROPOSITION 2.1. *For $\boldsymbol{\alpha} \in \mathbb{R}^k$ with $\lambda = \|\boldsymbol{\alpha}\| \neq 0$, for all $p > 2$,*

$$\sum_{i=1}^k \frac{|\alpha_i|^p}{\lambda^p} \leq 1,$$

*with equality if and only if $\boldsymbol{\alpha}$ is a multiple of a standard basis vector. In the case $p = 3$, yielding $\varphi$ of (7),*

$$\frac{1}{\sqrt{k}} \leq \varphi \leq 1, \tag{16}$$

*with equality to the upper bound if and only if $\boldsymbol{\alpha}$ is a multiple of a standard basis vector, and equality to the lower bound if and only if $|\alpha_i| = |\alpha_j|$ for all $i, j$. In addition, when $\alpha_i \geq 0$ with $\sum_{i=1}^n \alpha_i = 1$, then*

$$\lambda \leq \varphi, \tag{17}$$

*with equality if and only $\boldsymbol{\alpha}$ is equal to a standard basis vector.*

PROOF.   Since $|\alpha_i|/\lambda \leq 1$ we have $|\alpha_i|^{p-2}/\lambda^{p-2} \leq 1$, yielding

$$\sum_{i=1}^{k} \frac{|\alpha_i|^p}{\lambda^p} = \sum_{i=1}^{k} \left( \frac{|\alpha_i|^{p-2}}{\lambda^{p-2}} \right) \frac{|\alpha_i|^2}{\lambda^2} \leq \sum_{i=1}^{k} \frac{\alpha_i^2}{\lambda^2} = 1,$$

with equality if and only if for some $i$ we have $|\alpha_i| = \lambda$, and $\alpha_j = 0$ for all $j \neq i$. By Hölder's inequality with $p = 3, q = 3/2$, we have

$$\left( \sum_{i=1}^{k} 1 \cdot \alpha_i^2 \right)^{3/2} \leq \sqrt{k} \sum_{i=1}^{k} |\alpha_i|^3,$$

giving the lower bound (16), with equality if and only if $\alpha_i^2$ is proportional to 1 for all $i$. For the claim (17), by considering the inequality between the squared mean and variance of a random variable which takes the value $\alpha_i$ with probability $\alpha_i$, we have $(\sum_i \alpha_i^2)^2 \leq \sum_i \alpha_i^3$, with equality if and only if the variable is constant.   $\square$

Lemma 2.2 shows how zero biasing an independent sum behaves like a contraction mapping.

LEMMA 2.2.   *For $\boldsymbol{\alpha} \in \mathbb{R}^k$ with $\lambda = \|\boldsymbol{\alpha}\| \neq 0$, let*

$$Y = \sum_{i=1}^{k} \frac{\alpha_i}{\lambda} W_i,$$

*where $W_i$ are mean zero, variance 1, independent random variables distributed as $W$. Then*

$$d(Y, Y^*) \leq \varphi \, d(W, W^*)$$

*with $\varphi$ as in (7), and $\varphi < 1$ if and only if $\boldsymbol{\alpha}$ is not a multiple of a standard basis vector.*

PROOF.   By [1], for any collection $W_i^*$ with the $W_i$ zero biased distribution independent of $W_j$, $j \neq i$, and $I$ a random index independent of all other variables with distribution

$$P(I = i) = \frac{\alpha_i^2}{\lambda^2},$$

the variable

(18)                    $$Y^* = Y - \frac{\alpha_I}{\lambda}(W_I - W_I^*)$$

has the $Y$ zero biased distribution. Since $W_i \overset{d}{=} W$, we may take $(W_i, W_i^*) \overset{d}{=} (W, W^*)$, with $W, W^*$ achieving the infimum in (15). Then

$$|Y - Y^*| = \sum_{i=1}^{k} \frac{|\alpha_i|}{\lambda} |W_i - W_i^*| \mathbf{1}(I = i)$$

and

$$E|Y - Y^*| = \sum_{i=1}^{k} \frac{|\alpha_i|^3}{\lambda^3} E|W_i - W_i^*| = \left( \sum_{i=1}^{k} \frac{|\alpha_i|^3}{\lambda^3} \right) E|W - W^*|.$$

Now using (15) to upper bound $d(Y, Y^*)$ by the particular coupling in (18) we obtain

$$d(Y, Y^*) \leq E|Y - Y^*| = \varphi E|W - W^*| = \varphi \, d(W, W^*).$$

The final claim was shown in Proposition 2.1. $\quad\square$

In the classical case, when $Y = n^{-1/2} \sum_{i=1}^{n} W_i$, the normalized sum of i.i.d. random variables, applying Lemmas 2.1 and 2.2 with $\alpha_i = 1/\sqrt{n}$ gives $d(Y, \mathcal{N}) \leq 2d(Y, Y^*) \leq 2n^{-1/2} d(W, W^*) \to 0$ as $n \to \infty$, yielding a streamlined proof of the central limit theorem, complete with a bound in $d$.

When the sequence $X_n$ is given by the recursion (6) with $Z_n = 0$, setting $\lambda_n = \|\boldsymbol{\alpha}_n\|$ and $\sigma_n^2 = \mathrm{Var}(X_n)$ we have $\sigma_{n+1} = \lambda_n \sigma_n$, and we can write (6) as

$$W_{n+1} = \sum_{i=1}^{k} \frac{\alpha_{n,i}}{\lambda_n} W_{n,i} \qquad \text{with } W_n = \frac{X_n - c_n}{\sigma_n}.$$

Iterating the bound provided by Lemma 2.2 gives

$$d(W_n, W_n^*) \leq \left( \prod_{i=0}^{n-1} \varphi_i \right) d(W_0, W_0^*),$$

where

(19)
$$\varphi_n = \left( \frac{\sum_{i=1}^{k} |\alpha_{i,n}|^3}{\lambda_n^3} \right).$$

When $\limsup_n \varphi_n = \varphi < 1$, for any $\gamma \in (\varphi, 1)$ there exists $C$ such that for all $n$ we have $d(W_n, \mathcal{N}) \leq 2d(W_n, W_n^*) \leq C\gamma^n$. In Section 3 we study the situation when $Z_n$ is not necessarily zero.

**3. Bounds to the normal for approximately linear recursions.** In this section we study sequences $\{X_n\}_{n \geq 0}$ generated by the approximate linear recursion (6), and we present Theorem 3.1, which shows the exponential bound (5) holds when the perturbation term $Z_n$ is small as reflected in the term $\beta_n$ of (24), and holds in particular under the moment bounds in Conditions 3.1 and 3.2. When $Z_n$ is small, $X_{n+1}$ will be approximately equal to $\boldsymbol{\alpha}_n \cdot \mathbf{X}_n$, and therefore its variance $\sigma_{n+1}^2$ will be close to $\sigma_n^2 \lambda_n^2$, where $\lambda_n = \|\boldsymbol{\alpha}_n\|$, and the ratio $(\lambda_n \sigma_n)/\sigma_{n+1}$ will be close to 1. Iterating, the variance of $X_n$ will grow like a constant $C$ times $\lambda_{n-1}^2 \cdots \lambda_0^2$, so when $c_n \to c$ and $\boldsymbol{\alpha}_n \to \boldsymbol{\alpha}$, like $C^2 \lambda^{2n}$. Condition 3.1 assures that $Z_n$ is small relative to $X_n$ in that its variance grows at a slower rate. This condition was assumed in [13] for deriving a normal limiting law for the standardized sequence generated by (6).

CONDITION 3.1. The nonzero sequence of vectors $\boldsymbol{\alpha}_n \in \mathbb{R}^k, k \geq 2$, converges to $\boldsymbol{\alpha}$, not equal to any multiple of a standard basis vector. For $\lambda = \|\boldsymbol{\alpha}\|$, there exist $0 < \delta_1 < \delta_2 < 1$ and constants $C_{Z,2}, C_{X,2}$ such that, for all $n$,

$$\mathrm{Var}(Z_n) \leq C_{Z,2}^2 \lambda^{2n} (1 - \delta_2)^{2n},$$

$$\mathrm{Var}(X_n) \geq C_{X,2}^2 \lambda^{2n} (1 - \delta_1)^{2n}.$$

Bounds on the distance between $X_n$ and the normal can be provided under the following conditions on the fourth-order moments of $X_n$ and $Z_n$.

CONDITION 3.2. There exist $\delta_3$ and $\delta_4 \in (\delta_1, 1)$ and constants $C_{Z,4}, C_{X,4}$ such that

$$E(Z_n - EZ_n)^4 \leq C_{Z,4}^4 \lambda^{4n} (1 - \delta_4)^{4n},$$

$$E(X_n - EX_n)^4 \leq C_{X,4}^4 \lambda^{4n} (1 + \delta_3)^{4n}$$

and

(20)
$$\beta = \max\{\phi_1, \phi_2\} < 1$$

where $\phi_1 = \dfrac{(1 - \delta_2)(1 + \delta_3)^3}{(1 - \delta_1)^4}$ and $\phi_2 = \left(\dfrac{1 - \delta_4}{1 - \delta_1}\right)^2$.

Using Hölder's inequality and Condition 3.2 we may take

(21)
$$1 - \delta_2 \leq 1 - \delta_4 < 1 - \delta_1 \leq 1 + \delta_3.$$

In particular, $\beta \leq \eta$ for

(22)
$$\eta = \frac{(1 - \delta_4)(1 + \delta_3)^3}{(1 - \delta_1)^4}.$$

THEOREM 3.1. Let $X_{n+1} = \boldsymbol{\alpha}_n \cdot \mathbf{X}_n + Z_n$ with $\lambda_n = \|\boldsymbol{\alpha}_n\| \neq 0$ and $\mathbf{X}_n$ a vector in $\mathbb{R}^k$ with i.i.d. components distributed as $X_n$ with mean $c_n$ and nonzero variance $\sigma_n^2$. Set

(23)
$$W_n = \frac{X_n - c_n}{\sigma_n}, \qquad Y_n = \sum_{i=1}^k \frac{\alpha_{n,i}}{\lambda_n} W_{n,i}$$

*and*

$$\beta_n = E|W_{n+1} - Y_n| + \tfrac{1}{2} E|W_{n+1}^3 - Y_n^3|. \tag{24}$$

*If there exist $(\beta, \varphi) \in (0,1)^2$ such that*

$$\limsup_{n \to \infty} \frac{\beta_n}{\beta^n} < \infty \tag{25}$$

*and $\varphi_n$ in (19) satisfies*

$$\limsup_{n \to \infty} \varphi_n = \varphi, \tag{26}$$

*then with $\gamma = \beta$ when $\varphi < \beta$, and for any $\gamma \in (\varphi, 1)$ when $\beta \le \varphi$, there exists $C$ such that*

$$d(W_n, \mathcal{N}) \le C\gamma^n. \tag{27}$$

*Under Conditions 3.1 and 3.2, (27) holds for all $\gamma \in (\max(\beta, \varphi), 1)$, with $\beta$ as in (20), and $\varphi = \sum_{i=1}^k |\alpha_i|^3 / \lambda^3 < 1$.*

PROOF. By Lemma 2.1, it suffices to prove the bound (27) holds for $d(W_n, W_n^*)$. Let $f \in \mathcal{F}$ with $\mathcal{F}$ given by (13). Then $|f''(x)| \le 1, |f'(x)| \le |x|, |f(x)| \le x^2/2$, and for $h$ given by (8) with $\sigma^2 = 1$ and the chosen $f$, differentiation of (8) yields

$$h'(w) = f''(w) - wf'(w) - f(w),$$

and therefore

$$|h'(w)| \le (1 + \tfrac{3}{2}w^2). \tag{28}$$

Letting $r_n = (\lambda_n \sigma_n)/\sigma_{n+1}$ and using (23), write $X_{n+1} = \boldsymbol{\alpha}_n \cdot \mathbf{X}_n + Z_n$ as

$$W_{n+1} = r_n Y_n + T_n, \qquad \text{where } T_n = \frac{\sigma_n}{\sigma_{n+1}} \left( \frac{Z_n - EZ_n}{\sigma_n} \right). \tag{29}$$

Now by (28) and the definition of $\beta_n$ in (24),

$$E|h(W_{n+1}) - h(Y_n)| = E \left| \int_{Y_n}^{W_{n+1}} h'(u)\, du \right| \le \beta_n.$$

From (10) with $\sigma^2 = 1$, using $\text{Var}(Y_n) = 1$,

$$
\begin{aligned}
|Ef'(W_{n+1}) - Ef'(W_{n+1}^*)| &= |Eh(W_{n+1}) - Nh| \\
&= |E(h(W_{n+1}) - h(Y_n) + h(Y_n) - Nh)| \\
&\le \beta_n + |Eh(Y_n) - Nh| \\
&\le \beta_n + |E(f'(Y_n) - f'(Y_n^*))| \\
&\le \beta_n + d(Y_n, Y_n^*) \qquad \text{[by (14)]} \\
&\le \beta_n + \varphi_n d(W_n, W_n^*) \qquad \text{(by Lemma 2.2).}
\end{aligned}
$$

Taking the supremum over $f \in \mathcal{F}$ on the left-hand side, using (14) again and letting $d_n = d(W_n, W_n^*)$ we obtain, for all $n \geq 0$,

$$d_{n+1} \leq \beta_n + \varphi_n d_n.$$

Iteration yields that, for all $n, n_0 \geq 0$,

$$(30) \qquad d_{n_0+n} \leq \sum_{j=n_0}^{n_0+n-1} \left( \prod_{i=j+1}^{n_0+n-1} \varphi_i \right) \beta_j + \left( \prod_{i=n_0}^{n_0+n-1} \varphi_i \right) d_{n_0}.$$

Now suppose the bounds (25) and (26) hold and recall the choice of $\gamma$. When $\varphi < \beta$ take $\bar{\varphi} \in (\varphi, \beta)$ so that $\varphi < \bar{\varphi} < \beta = \gamma$; when $\beta \leq \varphi$ set $\bar{\varphi} \in (\varphi, \gamma)$ so that $\beta \leq \varphi < \bar{\varphi} < \gamma$. Then for any $\overline{B}$ greater than the lim sup in (25) there exists $n_0$ such that, for all $n \geq n_0$,

$$\beta_n \leq \overline{B} \beta^n \quad \text{and} \quad \varphi_n \leq \bar{\varphi}.$$

Applying these inequalities in (30) and summing yields, for all $n \geq 0$,

$$d_{n+n_0} \leq \overline{B} \beta^{n_0} \left( \frac{\beta^n - \bar{\varphi}^n}{\beta - \bar{\varphi}} \right) + \bar{\varphi}^n d_{n_0};$$

since $\max(\beta, \bar{\varphi}) \leq \gamma$, (27) follows.

To prove the final claim it suffices to show that, under Conditions 3.1 and 3.2, (25) and (26) hold with $\beta < 1$ as defined in (20), and with $\varphi = \sum_{i=1}^{k} |\alpha_i|^3 / \lambda^3 < 1$. Lemma 6 of [13] gives that the limit as $n \to \infty$ of $\sigma_n / (\lambda_0 \cdots \lambda_{n-1})$ exists in $(0, \infty)$, and therefore

$$(31) \qquad \lim_{n \to \infty} r_n = 1 \quad \text{and} \quad \lim_{n \to \infty} \frac{\sigma_{n+1}}{\sigma_n} = \lambda.$$

Referring to the definition of $T_n$ in (29) and using (31) and Conditions 3.1 and 3.2, there exist $C_{t,2}, C_{t,4}$ such that

$$(E|T_n|)^2 \leq ET_n^2 = \mathrm{Var}(T_n) = \left( \frac{\sigma_n}{\sigma_{n+1}} \right)^2 \frac{\mathrm{Var}(Z_n)}{\mathrm{Var}(X_n)} \leq C_{t,2}^2 \left( \frac{1 - \delta_2}{1 - \delta_1} \right)^{2n},$$

$$ET_n^4 = \left( \frac{\sigma_n}{\sigma_{n+1}} \right)^4 E\left( \frac{Z_n - EZ_n}{\sigma_n} \right)^4 \leq C_{t,4}^4 \left( \frac{1 - \delta_4}{1 - \delta_1} \right)^{4n}.$$

By independence, a simple bound and Condition 3.2 for the inequality,

$$(E|Y_n|)^2 \leq EY_n^2 = \mathrm{Var}(Y_n) = 1,$$

$$EY_n^4 \leq 6E\left( \frac{X_n - c_n}{\sigma_n} \right)^4 \leq 6C_{X,4}^4 \left( \frac{1 + \delta_3}{1 - \delta_1} \right)^{4n}.$$

From (6), with $\sigma_{Z_n} = \sqrt{\mathrm{Var}(Z_n)}$, $\sigma_{n+1} \leq \lambda_n \sigma_n + \sigma_{Z_n}$ and $\lambda_n \sigma_n \leq \sigma_{n+1} + \sigma_{Z_n}$; hence with $C_{r,1} = C_{t,2}$ we have

$$|\lambda_n \sigma_n - \sigma_{n+1}| \leq \sigma_{Z_n} \quad \text{so} \quad |r_n - 1| \leq C_{r,1} \left( \frac{1 - \delta_2}{1 - \delta_1} \right)^n.$$

Since $|r_n^p - 1| \leq \sum_{j \geq 1} \binom{p}{j} |r_n - 1|^j$, using (21) there are $C_{r,p}$ such that

$$|r_n^p - 1| \leq C_{r,p} \Big( \frac{1 - \delta_2}{1 - \delta_1} \Big)^n, \qquad p = 1, 2, \ldots .$$

Now considering the first term of $\beta_n$ of (24), recalling (29),

$$E|W_{n+1} - Y_n| = E|(r_n - 1)Y_n + T_n|$$

$$\leq |r_n - 1| E|Y_n| + E|T_n| \leq (C_{r,1} + C_{t,2}) \Big( \frac{1 - \delta_2}{1 - \delta_1} \Big)^n,$$

which is upper bounded by a constant times $\phi_1^n$.

For the second term of (24) we have

$$E|W_{n+1}^3 - Y_n^3| = E|(r_n^3 - 1)Y_n^3 + 3r_n^2 Y_n^2 T_n + 3r_n Y_n T_n^2 + T_n^3|.$$

Using the triangle inequality, the first term is bounded by a constant times $\phi_1^n$ as

$$|r_n^3 - 1| E|Y_n^3| \leq |r_n^3 - 1|(EY_n^4)^{3/4} \leq 6^{3/4} C_{r,3} C_{X,4}^3 \Big( \frac{(1 - \delta_2)(1 + \delta_3)^3}{(1 - \delta_1)^4} \Big)^n.$$

Since $r_n \to 1$ by (31), it suffices to bound the next two terms without the factor of $r_n$. Thus,

$$E|Y_n^2 T_n| \leq \sqrt{EY_n^4 ET_n^2} \leq 6^{1/2} C_{X,4}^2 C_{t,2} \Big( \frac{(1 - \delta_2)(1 + \delta_3)^2}{(1 - \delta_1)^3} \Big)^n,$$

which is less than a constant times $\phi_1^n$ by (21), and finally,

$$E|Y_n T_n^2| \leq \sqrt{EY_n^2 ET_n^4} \leq C_{t,4}^2 \Big( \frac{1 - \delta_4}{1 - \delta_1} \Big)^{2n} \leq C_{t,4}^2 \phi_2^n,$$

$$E|T_n^3| \leq (ET_n^4)^{3/4} \leq C_{t,4}^3 \Big( \frac{1 - \delta_4}{1 - \delta_1} \Big)^{3n} \leq C_{t,4}^3 \phi_2^{3n/2}.$$

Hence (25) holds with the given $\beta$.

Since $\boldsymbol{\alpha}_n \to \boldsymbol{\alpha}$, we have $\varphi_n \to \varphi$. Under Condition 3.1, $\boldsymbol{\alpha}$ is not a scalar multiple of a standard basis vector and $\varphi < 1$ by Lemma 2.2. We finish by invoking the first part of the theorem. $\quad\square$

**4. Normal bounds for hierarchical sequences.** The following result, extending Proposition 9 of [13] to higher orders, is used to show that the moment bounds of Conditions 3.1 and 3.2 are satisfied under the hypotheses of Theorem 1.1, so that Theorem 3.1 may be invoked. The dependence of the constants in (33) and (34) on $\varepsilon$ is suppressed for notational simplicity.

PROPOSITION 4.1. *Let the hypotheses of Theorem* 1.1 *hold. Following* (6), *with* $c_n = EX_n$ *and* $\boldsymbol{\alpha}_n = F'(c_n)$, *define*

$$(32) \qquad Z_n = F(\mathbf{X}_n) - \boldsymbol{\alpha}_n \cdot \mathbf{X}_n.$$

*Then with* $\boldsymbol{\alpha}$ *the limit of* $\boldsymbol{\alpha}_n$ *and* $\lambda = \|\boldsymbol{\alpha}\|$, *for any* $p \geq 1$ *and* $\varepsilon > 0$, *there exist constants* $C_{X,p}, C_{Z,p}$ *such that*

$$(33) \qquad E|Z_n - EZ_n|^p \leq C_{Z,p}^p (\lambda + \varepsilon)^{2pn} \qquad \text{for all } n \geq 0$$

*and*

$$(34) \qquad E|X_n - c_n|^p \leq C_{X,p}^p (\lambda + \varepsilon)^{pn} \qquad \text{for all } n \geq 0.$$

PROOF. Expanding $F(\mathbf{X}_n)$ around $\mathbf{c}_n$, with $\boldsymbol{\alpha}_n = F'(\mathbf{c}_n)$,

$$(35) \qquad F(\mathbf{X}_n) = F(\mathbf{c}_n) + \sum_{i=1}^k \alpha_{n,i}(X_{n,i} - c_n) + R_2(\mathbf{c}_n, \mathbf{X}_n),$$

where

$$R_2(\mathbf{c}_n, \mathbf{X}_n) = \sum_{i,j=1}^k \int_0^1 (1-t) \frac{\partial^2 F}{\partial x_i \partial x_j}(\mathbf{c}_n + t(\mathbf{X}_n - \mathbf{c}_n))(X_{n,i} - c_n)(X_{n,j} - c_n)\, dt.$$

Since the second partials of $F$ are continuous on $\mathcal{D} = [a,b]^k$, with $\|\cdot\|$ the supremum norm on $\mathcal{D}$, $B = 2^{-1} \max_{i,j} \|\partial^2 F / \partial x_i \partial x_j\| < \infty$, we have

$$(36) \qquad |R_2(\mathbf{c}_n, \mathbf{X}_n)| \leq B \sum_{i,j=1}^k |(X_{n,i} - c_n)(X_{n,j} - c_n)|.$$

Using (32), (35) and (36), we have, for all $p \geq 1$,

$$(37) \qquad \begin{aligned} &E|Z_n - EZ_n|^p \\ &= E\left|F(\mathbf{X}_n) - c_{n+1} - \sum_{i=1}^k \alpha_{n,i}(X_{n,i} - c_n)\right|^p \\ &= E|F(\mathbf{c}_n) - c_{n+1} + R_2(\mathbf{c}_n, \mathbf{X}_n)|^p \\ &\leq 2^{p-1}\left(|F(\mathbf{c}_n) - c_{n+1}|^p + B^p E\left(\sum_{i,j} |(X_{n,i} - c_n)(X_{n,j} - c_n)|\right)^p\right). \end{aligned}$$

For the first term of (37), again using (36),

$$\begin{aligned} |F(\mathbf{c}_n) - c_{n+1}|^p &= |F(\mathbf{c}_n) - EF(\mathbf{X}_n)|^p = |ER_2(\mathbf{c}_n, \mathbf{X}_n)|^p \\ &\leq B^p \left(E \sum_{i,j} |(X_{n,i} - c_n)(X_{n,j} - c_n)|\right)^p \end{aligned}$$

$$(38) \qquad \leq B^p k^p \left( E \sum_{i=1}^{k} (X_{n,i} - c_n)^2 \right)^p$$

$$\leq B^p k^{2p} [E(X_n - c_n)^2]^p$$

$$\leq B^p k^{2p} E(X_n - c_n)^{2p},$$

using Hölder's inequality for the final step.

Similarly, for the expectation in (37),

$$E \left( \sum_{i,j} |(X_{n,i} - c_n)(X_{n,j} - c_n)| \right)^p \leq k^p E \left( \sum_{i=1}^{k} (X_{n,i} - c_n)^2 \right)^p$$

$$(39) \qquad \leq k^{2p-1} E \left( \sum_{i=1}^{k} (X_{n,i} - c_n)^{2p} \right)$$

$$= k^{2p} E(X_n - c_n)^{2p}.$$

Applying (38) and (39) in (37) we obtain for all $p \geq 1$, with $C_p = 2^p B^p k^{2p}$,

$$(40) \qquad E|Z_n - EZ_n|^p \leq C_p E(X_n - c_n)^{2p}.$$

It therefore suffices to prove (34) to demonstrate the proposition.

In Lemma 8 of [13], it is shown that when $F : [a, b]^k \to [a, b]$ is an averaging function and there exists $c \in [a, b]$ such that $X_n \xrightarrow{p} c$, then

$$(41) \qquad \text{for every } \varepsilon \in (0, 1) \text{ there exists } M \text{ such that, for all } n,$$

$$P(|X_n - c| > \varepsilon) \leq M \varepsilon^n.$$

Hence the large deviation estimate (41) holds under the given assumptions, and so also with $c_n$ replacing $c$ when $c_n \to c$. Since $X_n \in [a, b]$ and $X_n \xrightarrow{p} c$, $c_n = EX_n \to c$ by the bounded convergence theorem.

We now show that if $a_n$, $n = 0, 1, \ldots$, is a sequence such that for every $\varepsilon > 0$ there exists $M$ such that, for all $n \geq n_0$,

$$(42) \qquad a_{n+1} \leq (\lambda + \varepsilon)^p a_n + M(\lambda + \varepsilon)^{p(n+1)},$$

then for all $\varepsilon > 0$ there exists $C$ such that

$$(43) \qquad a_n \leq C(\lambda + \varepsilon)^{pn} \qquad \text{for all } n.$$

Let $\varepsilon > 0$ be given, and let $M$ and $n_0$ be such that (42) holds with $\varepsilon$ replaced by $\varepsilon/2$. Setting

$$\rho = 1 - \left( \frac{\lambda + \varepsilon/2}{\lambda + \varepsilon} \right)^p \quad \text{and} \quad C = \max \left\{ \frac{a_{n_0}}{(\lambda + \varepsilon)^{n_0}}, \frac{M}{\rho} \left[ \frac{\lambda + \varepsilon/2}{\lambda + \varepsilon} \right]^{p(n_0+1)} \right\},$$

it is trivial that (43) holds for $n = n_0$. Since the second quantity in the maximum decreases when $n_0$ is replaced by $n \geq n_0$, induction shows (43) holds for all $n \geq n_0$. By increasing $C$ if necessary, we have that (43) holds for all $n$.

Unqualified statements in the remainder of the proof below involving $\varepsilon$ and $M$ are to be read to mean that for every $\varepsilon > 0$ there exists $M$, not necessarily the same at each occurrence, such that the statement holds for all $n$. By (41),

$$
\begin{aligned}
E(X_n - c_n)^{2p} &= E[(X_n - c_n)^{2p}; |X_n - c_n| \leq \varepsilon] \\
&\quad + E[(X_n - c_n)^{2p}; |X_n - c_n| > \varepsilon] \\
&\leq \varepsilon^p E|X_n - c_n|^p + M\varepsilon^n,
\end{aligned}
$$

so from (40),

$$
(44) \qquad E|Z_n - EZ_n|^p \leq \varepsilon E|X_n - c_n|^p + M\varepsilon^n.
$$

Since for all $w, z$,

$$
|w + z|^p \leq (1 + \varepsilon)|w|^p + M|z|^p,
$$

definition (32) yields

$$
(45) \quad E|X_{n+1} - c_{n+1}|^p \leq (1 + \varepsilon)E\left|\sum_{i=1}^{k} \alpha_{n,i}(X_n - c_n)\right|^p + ME|Z_n - EZ_n|^p.
$$

Specializing (45) to the case $p = 2$ gives, for all $n$ sufficiently large,

$$
E(X_{n+1} - c_{n+1})^2 \leq (\lambda + \varepsilon)^2 E(X_n - c_n)^2 + ME(Z_n - EZ_n)^2.
$$

Applying (44) with $p = 2$ to this inequality yields, for all $n$ sufficiently large,

$$
\begin{aligned}
E(X_{n+1} - c_{n+1})^2 &\leq (\lambda + \varepsilon)^2 E(X_n - c_n)^2 + M\varepsilon^{2n+2} \\
&\leq (\lambda + \varepsilon)^2 E(X_n - c_n)^2 + M(\lambda + \varepsilon)^{2(n+1)}.
\end{aligned}
$$

Hence, with $p = 2$, (42) and therefore (43) are true for $a_n = E(X_n - c_n)^2$, yielding (34) for $p = 2$. Now apply Hölder's inequality to prove the case $p = 1$.

Assume now that (34) is true for all $2 \leq q < p$ in order to induct on $p$. Expand the first term in (45), letting $\mathbf{p} = (p_1, \ldots, p_k)$ and $|\mathbf{p}| = \sum_i p_i$. Use the induction hypotheses, and Proposition 2.1 in (46), to obtain for all $n$ sufficiently large, with $A_{X,p} = \max_{q<p} C_{X,q}$ and $B_{X,p}^p = k^{p-1}A_{X,p}^p$,

$$
\begin{aligned}
E&\left|\sum_{i=1}^{k} \alpha_{n,i}(X_n - c_n)\right|^p \\
&\leq \sum_{i=1}^{k} |\alpha_{n,i}|^p E|X_{n,i} - c_n|^p + \sum_{|\mathbf{p}|=p, p_i < p} \binom{p}{\mathbf{p}} E\prod_{i=1}^{k} |\alpha_{n,i}|^{p_i}|X_{n,i} - c_n|^{p_i}
\end{aligned}
$$

$$\leq E|X_n - c_n|^p \sum_{i=1}^{k} |\alpha_{n,i}|^p + \sum_{|\mathbf{p}|=p, p_i < p} \binom{p}{\mathbf{p}} \prod_{i=1}^{k} |\alpha_{n,i}|^{p_i} C_{X,p_i}^{p_i} (\lambda + \varepsilon)^{p_i n}$$

$$(46) \qquad \leq E|X_n - c_n|^p \sum_{i=1}^{k} |\alpha_{n,i}|^p + A_{X,p}^p (\lambda + \varepsilon)^{pn} \sum_{|\mathbf{p}|=p} \binom{p}{\mathbf{p}} \prod_{i=1}^{k} |\alpha_{n,i}|^{p_i}$$

$$= E|X_n - c_n|^p \sum_{i=1}^{k} |\alpha_{n,i}|^p + A_{X,p}^p (\lambda + \varepsilon)^{pn} \left( \sum_{i=1}^{k} |\alpha_{n,i}| \right)^p$$

$$\leq \sum_{i=1}^{k} |\alpha_{n,i}|^p (E|X_n - c_n|^p + B_{X,p}^p (\lambda + \varepsilon)^{pn})$$

$$\leq (\lambda + \varepsilon)^p E|X_n - c_n|^p + B_{X,p}^p (\lambda + \varepsilon)^{p(n+1)}.$$

Applying (44) and (46) to (45) gives

$$E|X_{n+1} - c_{n+1}|^p \leq (\lambda + \varepsilon)^p E|X_n - c_n|^p + M(\lambda + \varepsilon)^{p(n+1)},$$

from which we can conclude (43) for $a_n = E|X_n - c_n|^p$, completing the induction on $p$. We conclude (34) holds for all $p \geq 1$. $\square$

PROOF OF THEOREM 1.1. By replacing $X_n$ and $F(\mathbf{x})$ by $X_n/F(\mathbf{1}_k)^n$ and $F(\mathbf{x})/F(\mathbf{1}_k)$, respectively, we may assume $F$ is averaging. By property 1 of averaging functions, $F(\mathbf{c}) = c$, and differentiation yields $\sum_{i=1}^{n} \alpha_i = 1$. By property 2, monotonicity, $\alpha_i \geq 0$, and (17) of Proposition 2.1 yields $0 < \lambda \leq \varphi < 1$.

Inspection of (22) shows that, for any $\eta \in (\lambda, 1)$, there exists $\delta_1$ and $\delta_3$ in $(0,1)$ and $\delta_4$ in $(\delta_1, 1 - \lambda)$ yielding $\eta$. For example, to achieve values arbitrarily close to $\lambda$ from above, take $\delta_1$ and $\delta_3$ close to zero and $\delta_4$ close to $1 - \lambda$ from below. Set $\delta_2 = \delta_4$. By Theorem 3.1 it suffices to show that Conditions 3.1 and 3.2 are satisfied for these choices of $\delta$.

Since $\delta_4 < 1 - \lambda$ we have $\lambda^2 < \lambda(1 - \delta_4)$; hence we may pick $\varepsilon > 0$ such that $(\lambda + \varepsilon)^2 < \lambda(1 - \delta_4)$. By Proposition 4.1, for $p = 2$ and $p = 4$, for this $\varepsilon$ there exists $C_{Z,p}^p$ such that

$$E(Z_n - EZ_n)^p \leq C_{Z,p}^p (\lambda + \varepsilon)^{4pn} \leq C_{Z,p}^p \lambda^{2pn} (1 - \delta_4)^{2pn}.$$

Hence the fourth and second moment bounds on $Z_n$ are satisfied with $\delta_4$ and $\delta_2 = \delta_4$, respectively.

Proposition 10 of [13] shows that under the assumptions of Theorem 1.1, for every $\varepsilon > 0$ there exists $C_{X,2}^2$ such that

$$\mathrm{Var}(X_n) \geq C_{X,2}^2 (\lambda - \varepsilon)^{2n}.$$

Taking $\varepsilon = \lambda \delta_1$, we have $\mathrm{Var}(X_n)$ satisfies its lower bound condition. Lastly, applying Proposition 4.1 with $p = 4$ and $\varepsilon = \lambda \delta_3$ we see the fourth moment bound on $X_n$ is satisfied, and the proof is complete. $\square$

**5. Convergence rates for the diamond lattice.** We now apply Theorem 1.1 to hierarchical sequences generated by the diamond lattice conductivity function $F$ in (2), for various choice of positive weights satisfying $F(\mathbf{1}_4) = 1$. For all such $F(\mathbf{x})$ the result of Shneiberg [8] quoted in Section 1 shows that $X_n$ satisfies a strong law if $X_0 \in [0, 1]$, say. The first partials of $F$ have the form, for example,

$$\frac{\partial F(\mathbf{x})}{\partial x_1} = \frac{(w_1 x_1^2)^{-1}}{((w_1 x_1)^{-1} + (w_2 x_2)^{-1})^2},$$

and therefore $F'(c_n \mathbf{1}_4)$ does not depend on $c_n$. In particular, for all $n$,

$$\boldsymbol{\alpha}_n = \left[ \frac{w_1^{-1}}{(w_1^{-1} + w_2^{-1})^2}, \frac{w_2^{-1}}{(w_1^{-1} + w_2^{-1})^2}, \frac{w_3^{-1}}{(w_3^{-1} + w_4^{-1})^2}, \frac{w_4^{-1}}{(w_3^{-1} + w_4^{-1})^2} \right]^\mathsf{T},$$

from which

$$(47) \qquad \varphi = \lambda^{-3} \left( \frac{w_1^{-3} + w_2^{-3}}{(w_1^{-1} + w_2^{-1})^6} + \frac{w_3^{-3} + w_4^{-3}}{(w_3^{-1} + w_4^{-1})^6} \right),$$

where

$$\lambda = \left( \frac{w_1^{-2} + w_2^{-2}}{(w_1^{-1} + w_2^{-1})^4} + \frac{w_3^{-2} + w_4^{-2}}{(w_3^{-1} + w_4^{-1})^4} \right)^{1/2}.$$

As an illustration, define the "side equally weighted network" to be the one with $\mathbf{w} = (w, w, 2 - w, 2 - w)^\mathsf{T}$ for $w \in [1, 2)$; such weights are positive and satisfy $F(\mathbf{1}_4) = 1$. For $w = 1$ all weights are equal, and we have $\boldsymbol{\alpha} = 4^{-1} \mathbf{1}_4$, and hence $\varphi$ achieves its minimum value $1/2 = 1/\sqrt{k}$ with $k = 4$. By Theorem 1.1, for all $\gamma \in (0, 1/2)$ there exists a constant $C$ such that $d(W_n, \mathcal{N}) \leq C\gamma^n$, with $\gamma$ close to $1/2$ corresponding to the rate $N^{-1/2+\varepsilon}$ for small $\varepsilon > 0$ and $N = 4^n$, the number of variables at stage $n$. As $w$ increases from 1 to 2, $\varphi$ increases continuously from $1/2$ to $1/\sqrt{2}$, with $w$ close to 2 corresponding to the least favorable rate for the side equally weighted network of $N^{-1/4+\varepsilon}$ for any $\varepsilon > 0$.

With only the restriction that the weights are positive and satisfy $F(\mathbf{1}_4) = 1$ consider

$$\mathbf{w} = (1 + 1/t, s, t, 1/t)^\mathsf{T},$$

$$\text{where } s = [(1 - (1/t + t)^{-1})^{-1} - (1 + 1/t)^{-1}]^{-1}, \ t > 0.$$

When $t = 1$ we have $s = 2/3$ and $\varphi = 11\sqrt{2}/27$. As $t \to \infty$, $s/t \to 1/2$ and $\boldsymbol{\alpha}$ tends to the standard basis vector $(1, 0, 0, 0)$, so $\varphi \to 1$. Since $11\sqrt{2}/27 < 1/\sqrt{2}$, the above two examples show that the value of $\gamma$ given by Theorem 1.1 for the diamond lattice can take any value in the range $(1/2, 1)$, corresponding to $N^{-\theta}$ for any $\theta \in (0, 1/2)$.

**6. Composition of strict averaging functions.** In this section, we prove Theorem 1.2, which shows when the composition of strictly averaging functions is again strictly averaging.

PROOF OF THEOREM 1.2. We first show $F_\mathbf{s}(\mathbf{x})$ satisfies the strict form of property 1. If $\mathbf{x} = t\mathbf{1}_k$, then $F_\mathbf{s}(\mathbf{x}) = F_0(s_1 t, \ldots, s_k t) = F_0(\mathbf{s})t = t$ and property 1 is satisfied in this case. Hence assume $\min_i x_i = x < y = \max_i x_i$. For such $\mathbf{x}$, if there is a $t$ such that $F_i(\mathbf{x}_i) = t$ for all $i = 1, \ldots, k$, then for some $i$ and $j$ we have $y = x_j$, $j \in I_i$, and hence $x < F_i(\mathbf{x}_i) = t$ since $F_i$ is strictly averaging, and similarly, $t < y$. Hence $x < F_\mathbf{1}(\mathbf{x}) = t < y$.

For $\mathbf{x}$ such that for all $i \in I_0$, $s_i F_i(\mathbf{x}_i) = t$ for some $t$, we have

$$F_\mathbf{s}(\mathbf{x}) = F_0(s_1 F_1(\mathbf{x}_1), \ldots, s_k F_k(\mathbf{x}_k)) = F_0(t\mathbf{1}_k) = t.$$

For $\mathbf{s} = \mathbf{1}_k$ we have just shown the strict inequality $x < t < y$ holds. Otherwise $\mathbf{s} \neq \mathbf{1}_k$ and by $F_0(\mathbf{s}) = 1$ we have $\min_i s_i < 1 < \max_i s_i$, and since $t = F_i(\mathbf{x}_i)/s_i$ for all $i$ there exist $i_1$ and $i_2$ such that

$$x \leq F_{i_1}(\mathbf{x}_{i_1}) < t < F_{i_2}(\mathbf{x}_{i_2}) \leq y,$$

yielding again the required strict inequality.

For $\mathbf{x}$ such that there are $i_1, i_2$ such that $s_{i_1} F_{i_1}(\mathbf{x}_{i_1}) \neq s_{i_2} F_{i_2}(\mathbf{x}_{i_2})$, we have $s_j F_j(\mathbf{x}_j) < \max_i s_i F_i(\mathbf{x}_i)$ for some $j$. Since $F_0$ is strictly monotone and homogeneous,

$$F_\mathbf{s}(\mathbf{x}) = F_0(s_1 F_1(\mathbf{x}_1), \ldots, s_k F_k(\mathbf{x}_k)) < F_0\left(s_1 \max_i F_i(\mathbf{x}_i), \ldots, s_k \max_i F_i(\mathbf{x}_i)\right)$$

$$= \max_i F_i(\mathbf{x}_i) F_0(\mathbf{s}) = \max_i F_i(\mathbf{x}_i) \leq y.$$

The argument for the minimum is the same; hence $F_\mathbf{s}(\mathbf{x})$ satisfies the strict form of property 1.

Since the composition of strictly monotone increasing functions is strictly monotone, the strict form of property 2 is satisfied for $F_\mathbf{s}(\mathbf{x})$.

The claim for $F_\mathbf{1}(\mathbf{x})$ now follows by setting

$$G_i(\mathbf{x}) = \frac{F_i(\mathbf{x}_i)}{F_i(\mathbf{1}_{|I_i|})} \quad \text{and} \quad s_i = \frac{F_i(\mathbf{1}_{|I_i|})F_0(\mathbf{1}_k)}{F_0(F_1(\mathbf{1}_{|I_1|}), \ldots, F_k(\mathbf{1}_{|I_k|}))}$$

$$\text{for } i = 0, 1, \ldots, k,$$

so that

$$\frac{F_\mathbf{1}(\mathbf{x})}{F_\mathbf{1}(\mathbf{1}_k)} = \frac{F_0(F_1(\mathbf{x}_1), \ldots, F_k(\mathbf{x}_k))}{F_0(F_1(\mathbf{1}_{|I_1|}), \ldots, F_k(\mathbf{1}_{|I_k|}))} = G_0(s_1 G_1(\mathbf{x}_1), \ldots, s_k G_k(\mathbf{x}_k)),$$

where $G_i(\mathbf{x}_i)$ is strictly averaging with $G_0$ homogeneous, and $G_0(\mathbf{s}) = 1$. □

DEPARTMENT OF MATHEMATICS
UNIVERSITY OF SOUTHERN CALIFORNIA
KAP 108
LOS ANGELES, CALIFORNIA 90089-2532
USA
E-MAIL: larry@math.usc.edu